\documentclass[12pt,reqno]{amsart}

\usepackage{amsmath}
\usepackage{amssymb}
\usepackage{amsthm}
\usepackage{graphicx}
\usepackage{geometry}
\usepackage{microtype}
\usepackage{booktabs}
\usepackage[hidelinks]{hyperref}

\newtheorem{theorem}{Theorem}[section]
\newtheorem{proposition}[theorem]{Proposition}
\newtheorem{lemma}[theorem]{Lemma}
\newtheorem{corollary}[theorem]{Corollary}
\newtheorem*{theoremA}{Theorem A}
\newtheorem*{theoremB}{Theorem B}

\theoremstyle{definition}
\newtheorem{definition}[theorem]{Definition}
\newtheorem{example}[theorem]{Example}
\newtheorem{remark}[theorem]{Remark}

\numberwithin{equation}{section}

\newcommand{\Q}{\mathbb{Q}}
\newcommand{\R}{\mathbb{R}}
\newcommand{\Z}{\mathbb{Z}}
\newcommand{\cD}{\mathcal{D}}
\newcommand{\cR}{\mathcal{R}}
\newcommand{\cF}{\mathcal{F}}
\newcommand{\Zt}{\Z_{(2)}}

\begin{document}

\title[Area obstructions for perfect triangles]
{Signed Admissible Regions and 2-adic Area\\
Obstructions for Perfect Triangles}

\author{Qizhen Wang}
% TODO: add the department and your e-mail address.
\address{Waseda University, Tokyo, Japan}
\email{}

\subjclass[2020]{Primary 11D09; Secondary 11D88, 51M04}

\keywords{perfect triangles, MED triangles, rational medians,
Heron triangles, Buchholz parametrization, signed admissibility,
2-adic valuations}

\begin{abstract}
A perfect triangle is a triangle whose sides, medians and area are
all integers; whether one exists is an open problem. We study
Buchholz's rational parametrization of triangles with two rational
medians on the whole $(\theta,\phi)$-plane, allowing the generated
side lengths to be negative. We factor its Heron polynomial
$H(\theta,\phi)$ into ten explicit factors and show that the points
whose componentwise absolute values form a triangle are exactly
those with $H>0$; the points giving positive side lengths directly
form three open quadrilaterals, one of which is the region used by
Buchholz and Rathbun. Every triangle with two rational medians
arises from sixteen signed parameter points, given by explicit
rational functions. Rational area forces $H$ to be a square in
$\Q_2$, and the valuation and the unit part of $H$ yield parity and
modulo-$8$ conditions in three valuation ranges. Since a triangle is
excluded as soon as one of its sixteen parameter points fails a
test, these conditions exclude infinite families of signed
admissible parameters, including non-integral ones. The conditions
are necessary only and do not settle the existence of perfect
triangles.
\end{abstract}

\maketitle

% =====================================================================
\section{Introduction}
\label{sec:introduction}
% =====================================================================

A \emph{MED triangle} is a non-degenerate triangle whose three sides
and three medians have integer lengths. A \emph{perfect triangle} is
a MED triangle whose area is also an integer; after scaling, this is
the same as a Heron triangle (rational sides and area) with three
rational medians. Whether a perfect triangle exists is an open
problem; see Guy~\cite[D21]{Guy2004}.

For a triangle with sides $a,b,c$ and medians $m_a,m_b,m_c$ we work
with the \emph{doubled medians}
\begin{equation}\label{eq:doubled-medians}
K^2=2b^2+2c^2-a^2,\qquad
L^2=2a^2+2c^2-b^2,\qquad
M^2=2a^2+2b^2-c^2,
\end{equation}
so that $K=2m_a$, $L=2m_b$ and $M=2m_c$. Integer medians force all
three sides to be even. For questions of existence, rational lengths
suffice: a triangle with rational sides and rational medians can be
scaled to a MED triangle, and if its area is also rational, a
suitable integer scaling makes the six lengths and the area integral.
Conversely, an irrational area stays irrational under every nonzero
rational scaling, since the area is multiplied by a rational square.

\subsection*{Previous work}
Euler~\cite{Euler} gave a parametric family of triangles with
rational sides and medians; not every such triangle belongs to
it~\cite{Buchholz2002}. Buchholz~\cite{Buchholz1989} found a rational
parametrization of \emph{all} triangles with two rational medians, in
two parameters $(\theta,\phi)$, and showed~\cite{Buchholz2002} that
imposing the third median leads to a one-parameter family of elliptic
curves. Buchholz and Rathbun~\cite{BR1997,BR1998} used the same
parametrization on the region $0<\theta,\phi<1$, $\phi+2\theta>1$
to construct infinitely many Heron triangles (rational sides and
area) with two rational medians. B\'acsk\'ai, Buchholz, Rathbun and
Smith~\cite{BBRS2003} studied the discrete symmetries of the problem,
including sign changes of the sides and medians and the exchange
$a\leftrightarrow b$, $k\leftrightarrow\ell$; in this way the main
family of Heron triangles gives rational points on eight isomorphic
genus-one curves $C_1,\dots,C_8$ in the $(\theta,\phi)$-plane
(see~\cite[Section~1]{Hone2022}). Buchholz and
Stingley~\cite{BS2019} showed that none of the triangles coming from
these curves has three rational medians, while four of the nine
smallest known examples do not lie on them~\cite{Hone2022}. Ismail~\cite{Ismail2020}
studied the perfect-triangle question on the curve $C_4$, and
Hone~\cite{Hone2022,Hone2024} related the Buchholz--Rathbun family
to Somos-$5$ sequences and surveyed the subject. The obstructions
in the present paper do not use these curves. Non-degenerate right and isosceles
triangles cannot be MED triangles; see
Buchholz~\cite{Buchholz1989}, Mordell~\cite[pp.~19--21]{Mordell1969}
and Dickson~\cite[p.~475]{Dickson1920}.

\subsection*{This paper}
Searches and constructions in the works cited above take place in
the region where all three generated side lengths are positive, while
sign changes enter through the symmetry group of~\cite{BBRS2003}. Our
starting point is to treat signed parameter points systematically, as
points of the plane in their own right. All the quantities in the
problem depend only on $a^2,b^2,c^2$: the medians
by~\eqref{eq:doubled-medians} and the area by Heron's formula. Hence a
parametrized triple with \emph{mixed signs}, that is, with some
generated side lengths negative, is just as good as the triple of its
absolute values. We therefore study Buchholz's parametrization on the
whole real plane, keeping track of two notions of admissibility
(Definition~\ref{def:admissible}):
\emph{direct} admissibility, where the parametrized entries themselves
form a triangle, and \emph{signed} admissibility, where their absolute
values do. Our first result describes both regions.

\begin{theoremA}
Let $H(\theta,\phi)$ be the Heron polynomial of Buchholz's
parametrization \eqref{eq:param}. Then
\[
H=16\,\phi(1-\phi)(1+\phi)\,\theta(\theta-1)(\theta+1)\,
G\,(\phi-\theta-1)(\theta+2\phi+1)(\phi+2\theta-1),
\]
where $G=(3\phi+1)\theta+1-\phi$, and:
\begin{enumerate}
\item the signed admissible region is $\cR=\{H>0\}$
  \textup{(}Proposition~\ref{prop:signed-region}\textup{)};
\item the direct admissible region is the union of three open
  quadrilaterals, one of which is the Buchholz--Rathbun region
  $\cF$ \textup{(}Proposition~\ref{prop:direct-region}\textup{)};
\item every triangle with two rational medians arises from exactly
  the sixteen points of $\cR\cap\Q^2$ described in
  Proposition~\ref{prop:sixteen}.
\end{enumerate}
\end{theoremA}

The same factorization gives arithmetic information. Rational area
requires $H$ to be a rational square, hence a square in $\R$ (this
is signed admissibility) and in $\Q_2$. Writing
$\mu=v_2(\theta)$, $\nu=v_2(\phi)$, $\kappa=\min\{\mu,\nu\}$ and
$\lambda=\max\{\mu,\nu\}$, we obtain the following.

\begin{theoremB}
Let $(\theta,\phi)\in\cR\cap\Q^2$. Then
\[
v_2(H)=
\begin{cases}
4+\mu+\nu, & \mu,\nu\ge1,\\
5+7\kappa+\lambda, & \kappa\le-2,\ \lambda\ge1,\\
5+7\kappa+4\lambda, & \lambda\le-1,\ \kappa\le\lambda-2.
\end{cases}
\]
Consequently the generated triangle has irrational area if
$\mu+\nu$ is odd in the first range, if $\mu+\nu$ is even in the
second, and if $\kappa$ is even in the third. Further conditions
modulo $8$ hold in the ranges $\mu,\nu\ge2$ and $\kappa\le-3$,
$\lambda\ge3$ \textup{(}Propositions~\ref{prop:valuation-obstructions},
\ref{prop:positive-unit} and~\ref{prop:mixed-unit}\textup{)}.
\end{theoremB}

A triangle is excluded as soon as \emph{one} of its sixteen parameter
points fails one of these tests (Corollary~\ref{cor:sixteen-tests}).
This yields explicit infinite families of signed admissible parameters
with irrational area (Corollary~\ref{cor:infinite-families}). The
mixed-sign points are essential here: for five of these families, in
the range we checked, none of the parameter points obtained from the
positive side lengths of the triangle is excluded by our tests, so the
exclusion is visible only through parameter points with mixed signs
(Remark~\ref{rem:signed-needed}). The conditions are necessary but not
sufficient, and they do not settle the existence of perfect
triangles.

Section~\ref{sec:parametrization} recalls the parametrization in a
form that can be checked by direct expansion, proves its
completeness, and discusses the third median and two examples.
Section~\ref{sec:regions} describes the admissible regions and the
sixteen parameter points of a triangle. Section~\ref{sec:area}
derives the $2$-adic obstructions. Buchholz's original derivation is
sketched in Appendix~\ref{app:derivation}.

% =====================================================================
\section{Buchholz's Parametrization}
\label{sec:parametrization}
% =====================================================================

\subsection{Signed triples}

It is convenient to allow the entries of a triple $(a,b,c)$ of real
numbers to have arbitrary signs. Its doubled medians are defined
by~\eqref{eq:doubled-medians}, and its \emph{Heron product} is
\begin{equation}\label{eq:heron-product}
\begin{aligned}
Q(a,b,c)&=(a+b+c)(-a+b+c)(a-b+c)(a+b-c)\\
&=2a^2b^2+2a^2c^2+2b^2c^2-a^4-b^4-c^4 .
\end{aligned}
\end{equation}
The expanded form shows $Q(a,b,c)=Q(|a|,|b|,|c|)$. If $|a|,|b|,|c|$
are the sides of a triangle with area $S$, then Heron's formula reads
$16S^2=Q(a,b,c)$.

\begin{lemma}\label{lem:heron-sign}
For real $a,b,c$, the numbers $|a|,|b|,|c|$ are the sides of a
non-degenerate triangle if and only if $Q(a,b,c)>0$.
\end{lemma}

\begin{proof}
By sign invariance we may assume $a,b,c\ge0$. If one entry is zero,
say $a=0$, then $Q=-(b^2-c^2)^2\le0$. If all entries are positive and
$c$ is the largest, then $a+b+c$, $-a+b+c$ and $a-b+c$ are positive,
so $Q$ has the sign of $a+b-c$. Hence $Q>0$ exactly when the strict
triangle inequalities hold.
\end{proof}

\subsection{The parametrization}

For $\theta,\phi\in\R$ put
\begin{equation}\label{eq:param}
\begin{aligned}
A_0&=-2\phi\theta^2+(-\phi^2+2\phi+1)\theta+1-\phi^2,\\
B_0&=(\phi-1)\theta^2+2(\phi^2+\phi)\theta+1-\phi,\\
C_0&=(\phi+1)\theta^2+(-\phi^2+2\phi+1)\theta+\phi^2-\phi,
\end{aligned}
\end{equation}
and, for a scaling factor $\tau\in\Q\setminus\{0\}$, consider the
triple $(a,b,c)=\tau(A_0,B_0,C_0)$. This is the parametrization of
Buchholz~\cite{Buchholz1989}, in the form used
in~\cite{BR1998,Hone2024}. We also put
\begin{equation}\label{eq:K0L0}
\begin{aligned}
K_0&=2\theta^2+(3\phi^2+2\phi+1)\theta-(\phi-1)^2,\\
L_0&=(3\phi+1)\theta^2+2(1-\phi)\theta+(1-\phi)(1+2\phi),
\end{aligned}
\end{equation}
and define the \emph{mixed factor}
\begin{equation}\label{eq:G}
G=(3\phi+1)\theta+1-\phi .
\end{equation}

\begin{proposition}\label{prop:param}
The following polynomial identities hold.
\begin{enumerate}
\item Two rational medians:
\[
2B_0^2+2C_0^2-A_0^2=K_0^2,\qquad
2A_0^2+2C_0^2-B_0^2=L_0^2 .
\]
\item Triangle factors:
\begin{equation}\label{eq:triangle-factors}
\begin{aligned}
T_1&:=A_0+B_0+C_0=2G,\\
T_2&:=-A_0+B_0+C_0=2\phi(\theta+1)(\phi+2\theta-1),\\
T_3&:=A_0-B_0+C_0=2\theta(1-\phi)(\theta+2\phi+1),\\
T_4&:=A_0+B_0-C_0=2(\theta-1)(\phi+1)(\phi-\theta-1).
\end{aligned}
\end{equation}
\item Inverse formulas: if $G\neq0$, then
\begin{equation}\label{eq:inverse-K0L0}
\theta=\frac{C_0-A_0+L_0}{A_0+B_0+C_0},\qquad
\phi=\frac{B_0-C_0+K_0}{A_0+B_0+C_0}.
\end{equation}
\end{enumerate}
In particular, for all $\theta,\phi,\tau\in\Q$ with $\tau\neq0$, the
triple $\tau(A_0,B_0,C_0)$ has rational doubled medians, namely
$K=\pm\tau K_0$ and $L=\pm\tau L_0$.
\end{proposition}

\begin{proof}
Each identity follows by expanding both sides.
\end{proof}

Proposition~\ref{prop:param} can be checked without any reference to
how the parametrization was found; Buchholz's derivation by norm
equations in $\Q(\sqrt2)$ is sketched in Appendix~\ref{app:derivation}.
Throughout, $\phi$ is the parameter attached to the median $K$ on the
first side and $\theta$ the parameter attached to the median $L$ on
the second side.

\subsection{Completeness}

We now show that every triple with two rational doubled medians is
proportional to a parametrized one. For a triple $(a,b,c)$ with
$s=a+b+c\neq0$, consider the quadratic polynomials
\begin{equation}\label{eq:quadratics}
\begin{aligned}
P_\theta(x)&=s\,x^2+2(a-c)\,x+(b-a-c)=n_1(x)\cdot(a,b,c),\\
P_\phi(y)&=s\,y^2-2(b-c)\,y+(a-b-c)=n_2(y)\cdot(a,b,c),
\end{aligned}
\end{equation}
where
\[
n_1(x)=(x^2+2x-1,\ x^2+1,\ x^2-2x-1),\qquad
n_2(y)=(y^2+1,\ y^2-2y-1,\ y^2+2y-1).
\]
Their quarter-discriminants are
\[
(a-c)^2-s(b-a-c)=L^2,\qquad (b-c)^2-s(a-b-c)=K^2,
\]
so their roots are
\begin{equation}\label{eq:inverse}
x=\frac{c-a\pm L}{a+b+c},\qquad
y=\frac{b-c\pm K}{a+b+c}.
\end{equation}
In particular the roots are rational if and only if $K$ and $L$ are
rational.

\begin{lemma}\label{lem:cross}
For all $\theta,\phi$,
\[
n_1(\theta)\times n_2(\phi)=-2\,\big(A_0,B_0,C_0\big)(\theta,\phi).
\]
Moreover, $A_0=B_0=C_0=0$ holds only at the three points
$(\theta,\phi)=(1,-1)$, $(0,1)$ and $(-1,0)$.
\end{lemma}

\begin{proof}
The first statement follows by expansion. If $A_0=B_0=C_0=0$, then
the four products in~\eqref{eq:triangle-factors} vanish. If
$\theta=1$ or $\phi=-1$, then $G=0$ gives $(\theta,\phi)=(1,-1)$. If
$\phi=\theta+1$, then $G=3\theta(\theta+1)$, so
$(\theta,\phi)=(0,1)$ or $(-1,0)$. Conversely, $A_0=B_0=C_0=0$ at
these three points.
\end{proof}

\begin{theorem}[Completeness]\label{thm:completeness}
Let $a,b,c\in\Q$ satisfy $Q(a,b,c)\neq0$, and suppose that $K$ and $L$
in~\eqref{eq:doubled-medians} are rational. Then $a+b+c\neq0$, the
four numbers in~\eqref{eq:inverse} are rational, and for each of the
four pairs $(\theta,\phi)$ formed by a root $\theta$ of $P_\theta$ and
a root $\phi$ of $P_\phi$ there is $\rho\in\Q\setminus\{0\}$ with
\[
\big(A_0,B_0,C_0\big)(\theta,\phi)=\rho\,(a,b,c).
\]
In particular, every triangle with rational sides and rational medians
on two sides is similar to the triangle with sides
$|A_0|,|B_0|,|C_0|$ for some $(\theta,\phi)\in\Q^2$.
\end{theorem}

\begin{proof}
The factor $a+b+c$ of $Q(a,b,c)$ is nonzero, and the roots are
rational by~\eqref{eq:inverse}. Fix a pair of roots $(\theta,\phi)$.
By~\eqref{eq:quadratics}, the nonzero vector $(a,b,c)$ is orthogonal
to $n_1(\theta)$ and $n_2(\phi)$. Suppose first that
$(A_0,B_0,C_0)(\theta,\phi)\neq0$. By Lemma~\ref{lem:cross},
$n_1(\theta)$ and $n_2(\phi)$ are then linearly independent, so
$(a,b,c)$ is a multiple of their cross product, that is, of
$(A_0,B_0,C_0)(\theta,\phi)$.

Otherwise $(\theta,\phi)$ is one of the three points in
Lemma~\ref{lem:cross}, and $\theta\in\{1,0,-1\}$. Then
$n_1(\theta)$ equals $(2,2,-2)$, $(-1,1,-1)$ or $(-2,2,2)$, so
$n_1(\theta)\cdot(a,b,c)=0$ gives $a+b-c=0$, $a-b+c=0$ or
$-a+b+c=0$. Each of these forces $Q(a,b,c)=0$, a contradiction.
\end{proof}

\begin{remark}
For a positive triangle, Theorem~\ref{thm:completeness} recovers
Buchholz's completeness statement. The proof uses only the signs of
$a,b,c$ through $Q$, so it applies verbatim to signed triples; this
is used in Section~\ref{subsec:sixteen}.
\end{remark}

\subsection{The third median}

Let $M_0^2=2A_0^2+2B_0^2-C_0^2$. Direct expansion gives
\begin{equation}\label{eq:quartic}
\begin{aligned}
M_0^2=F_\phi(\theta):={}&(3\phi-1)^2\theta^4
+2(9\phi^3-9\phi^2-11\phi-1)\theta^3\\
&+3(3\phi^4+6\phi^3+2\phi^2+2\phi-1)\theta^2\\
&+2(\phi-1)(3\phi^3-8\phi^2-11\phi-2)\theta
+(\phi-1)^2(\phi+2)^2 .
\end{aligned}
\end{equation}
Thus $\tau(A_0,B_0,C_0)$ has a rational third median if and only if
$F_\phi(\theta)$ is a rational square. The leading coefficient and the
constant term of $F_\phi$ are squares, so for fixed $\phi$ the curve
$\eta^2=F_\phi(\theta)$ has rational points at $\theta=0$ and at
infinity; when $F_\phi$ is squarefree it is an elliptic curve over
$\Q(\phi)$. This is the elliptic-curve viewpoint of
Buchholz~\cite{Buchholz2002}.

The third-median condition is independent of geometric admissibility
and of rational area. Sections~\ref{sec:regions} and~\ref{sec:area}
treat these two requirements without assuming it.

\subsection{Two MED triangles}
\label{subsec:med-examples}

\begin{figure}[t]
\centering
\includegraphics[width=\textwidth]{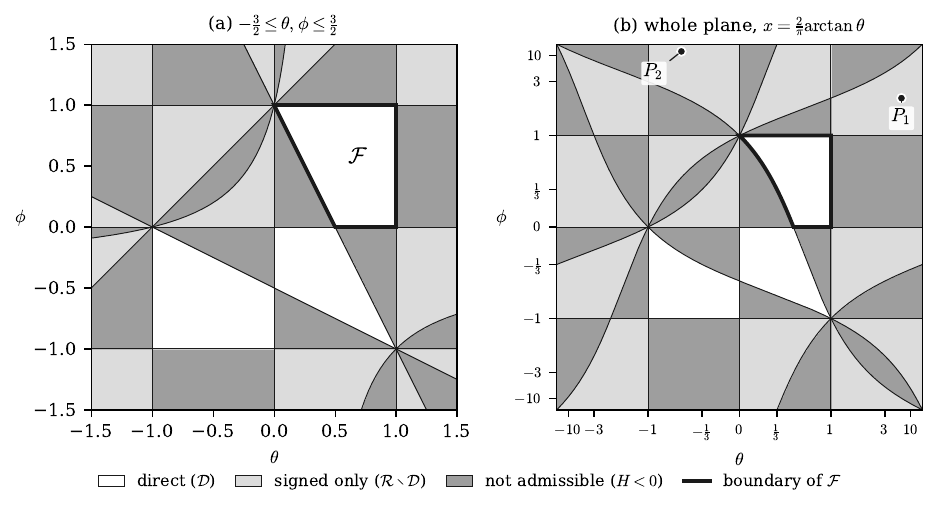}
\caption{Admissible regions. White: directly admissible ($\cD$).
Light gray: signed admissible only ($\cR\setminus\cD$). Dark gray:
$H<0$. Thin lines: the zero set of $H$. The thick outline is the
Buchholz--Rathbun region $\cF$. Panel~(a) shows
$-\frac32\le\theta,\phi\le\frac32$; panel~(b) shows the whole plane
in the coordinates $\frac{2}{\pi}\arctan\theta$ and
$\frac{2}{\pi}\arctan\phi$, with the points $P_1$ and $P_2$ of
Examples~\ref{ex:first-med} and~\ref{ex:second-med}.}
\label{fig:regions}
\end{figure}

The following two MED triangles illustrate how parameter points with
mixed-sign outputs represent genuine triangles. They are the points
$P_1$ and $P_2$ in Figure~\ref{fig:regions}.

\begin{example}[Euler's triangle]\label{ex:first-med}
Take $P_1=(\theta,\phi)=(11/2,\,2)$. Then
\[
F_2(\theta)=25\theta^4+26\theta^3+321\theta^2-64\theta+16,
\qquad F_2\!\left(\tfrac{11}{2}\right)=\left(\tfrac{765}{4}\right)^2,
\]
and $(A_0,B_0,C_0)=\tfrac38(-316,254,262)$. With $\tau=8/3$ we obtain
the signed triple $(a,b,c)=(-316,254,262)$ with doubled medians
$(K,L,M)=(408,522,510)$, up to sign. Thus the triangle with sides
$316,254,262$ has medians $204,261,255$. It is the example given by
Euler~\cite{Euler}. Since $316<254+262$, the point $P_1$ is signed
admissible, although its entries have mixed signs.
\end{example}

\begin{example}\label{ex:second-med}
Take $P_2=(\theta,\phi)=(-6/11,\,16)$. Then
\[
F_{16}(\theta)=2209\theta^4+68766\theta^3+665181\theta^2
+301860\theta+72900,
\qquad
F_{16}\!\left(-\tfrac{6}{11}\right)=\left(\tfrac{37332}{121}\right)^2,
\]
and $(A_0,B_0,C_0)=\tfrac{153}{242}(-226,-486,580)$. With
$\tau=242/153$ we obtain $(a,b,c)=(-226,-486,580)$ and, up to sign,
$(K,L,M)=(1046,734,488)$. Thus the triangle with sides $226,486,580$
has medians $523,367,244$. Since $580<226+486$, the point $P_2$ is
signed admissible but not directly admissible.
% TODO: add a reference in which this triangle appears, if desired.
\end{example}

% =====================================================================
\section{Admissible Regions in the Parameter Plane}
\label{sec:regions}
% =====================================================================

\begin{definition}\label{def:admissible}
A point $(\theta,\phi)\in\R^2$ is \emph{directly admissible} if
$A_0,B_0,C_0$ satisfy the three strict triangle inequalities, and
\emph{signed admissible} if $|A_0|,|B_0|,|C_0|$ are the sides of a
non-degenerate triangle. We write $\cD$ and $\cR$ for the two regions.
\end{definition}

Signed admissibility does not depend on the scaling factor $\tau\ne0$;
direct admissibility refers to $\tau>0$. Adding the two triangle
inequalities that contain a given entry shows that the entries of a
directly admissible triple are positive, so $\cD\subset\cR$. A triple
with three negative entries is not in $\cD$, although $\tau<0$ turns
it into a positive one; the examples of
Section~\ref{subsec:med-examples} have mixed signs and cannot be made
positive by any choice of $\tau$.

\subsection{Signed admissibility}

Recall from~\eqref{eq:triangle-factors} the four \emph{triangle
factors} of the parametrized triple,
\[
T_1=A_0+B_0+C_0,\quad
T_2=-A_0+B_0+C_0,\quad
T_3=A_0-B_0+C_0,\quad
T_4=A_0+B_0-C_0,
\]
and from~\eqref{eq:G} the \emph{mixed factor}
$G=(3\phi+1)\theta+1-\phi$, which satisfies $T_1=2G$. We call
\[
H(\theta,\phi)=Q(A_0,B_0,C_0)=T_1T_2T_3T_4
\]
the \emph{Heron polynomial} of the parametrization. Each $T_i$ contains
a factor $2$, which produces the constant $16=2^4$. Substituting the
factorizations~\eqref{eq:triangle-factors} gives
\begin{equation}\label{eq:H}
H=16\,\phi(1-\phi)(1+\phi)\,\theta(\theta-1)(\theta+1)\,
G\,(\phi-\theta-1)(\theta+2\phi+1)(\phi+2\theta-1).
\end{equation}

\begin{proposition}[Signed admissibility]
\label{prop:signed-region}
The signed admissible region is
$\cR=\{(\theta,\phi)\in\R^2: H(\theta,\phi)>0\}$.
\end{proposition}

\begin{proof}
Apply Lemma~\ref{lem:heron-sign} to $(A_0,B_0,C_0)$.
\end{proof}

The zero set of $H$ consists of the nine lines
\begin{equation}\label{eq:nine-lines}
\begin{gathered}
\phi=-1,\quad\phi=0,\quad\phi=1,\quad
\theta=-1,\quad\theta=0,\quad\theta=1,\\
\phi=\theta+1,\quad
\phi=-\tfrac{1}{2}(\theta+1),\quad
\phi=1-2\theta,
\end{gathered}
\end{equation}
together with the hyperbola $G=0$, that is,
$\phi=(\theta+1)/(1-3\theta)$ for $\theta\neq1/3$. On each connected
component of the complement the sign of $H$ is constant, and $\cR$
is the union of the components on which it is positive. The points of
Examples~\ref{ex:first-med} and~\ref{ex:second-med} satisfy
$H(P_1)>0$ and $H(P_2)>0$.

\subsection{Direct admissibility}

\begin{proposition}[Direct admissibility]
\label{prop:direct-region}
The direct admissible region is
\[
\cD=\{(\theta,\phi)\in\R^2:\ T_2>0,\ T_3>0,\ T_4>0\}.
\]
It is the union of the three open quadrilaterals with vertices
\[
\begin{aligned}
\cF&:\ (0,1),\ (\tfrac12,0),\ (1,0),\ (1,1);\\
\cD_2&:\ (-1,-1),\ (0,-1),\ (0,-\tfrac12),\ (-1,0);\\
\cD_3&:\ (0,0),\ (\tfrac12,0),\ (1,-1),\ (0,-\tfrac12).
\end{aligned}
\]
In particular $\cD\subset(-1,1)^2$, and
\[
\cD\cap(0,1)^2=\cF=\{0<\theta<1,\ 0<\phi<1,\ \phi+2\theta>1\},
\]
which is the region used by Buchholz and
Rathbun~\cite{BR1998} and by Hone~\cite{Hone2024}.
\end{proposition}

\begin{proof}
The first statement is a restatement of the triangle inequalities.
The factors $T_2,T_3,T_4$ have constant signs on each cell of the
arrangement of the nine lines~\eqref{eq:nine-lines}, so it suffices
to test one point in each of the finitely many cells; doing so
gives the three quadrilaterals. For the last statement, let
$0<\theta,\phi<1$. Then $\theta-1<0$, $\phi+1>0$ and $\phi-\theta-1<0$,
so $T_4>0$; also $\theta>0$, $1-\phi>0$ and $\theta+2\phi+1>0$, so
$T_3>0$; and $T_2>0$ if and only if $\phi+2\theta>1$.
\end{proof}

\subsection{The sixteen parameter points of a triangle}
\label{subsec:sixteen}

By Buchholz and Rathbun~\cite{BR1998}, the region $\cF$ already
parametrizes every triangle with two rational medians, so one may ask
what the rest of $\cR$ adds. The answer is that a triangle has
sixteen parameter points in $\cR$, and the arithmetic tests of
Section~\ref{sec:area} may be applied at any of them. The sign changes
and root choices below are part of the discrete symmetry group studied
by B\'acsk\'ai, Buchholz, Rathbun and Smith~\cite{BBRS2003}
(see also~\cite[Theorem~3.7]{Hone2022}); what we add is the
description of these points as rational functions on $\cR$, the fact
that they all lie in $\cR$, and the converse statement~(iv).

\begin{proposition}\label{prop:sixteen}
Let $(\theta,\phi)\in\cR\cap\Q^2$, and write $A_0,B_0,C_0,K_0,L_0$ for
the values at $(\theta,\phi)$. For
$\varepsilon=(\varepsilon_1,\varepsilon_2)\in\{\pm1\}^2$ put
\[
(a_\varepsilon,b_\varepsilon,c_\varepsilon)
=(\varepsilon_1A_0,\ \varepsilon_2B_0,\ C_0),
\qquad
s_\varepsilon=a_\varepsilon+b_\varepsilon+c_\varepsilon,
\]
and for $\delta=(\delta_1,\delta_2)\in\{\pm1\}^2$ put
\begin{equation}\label{eq:sixteen}
\theta_{\varepsilon,\delta}
=\frac{c_\varepsilon-a_\varepsilon+\delta_1L_0}{s_\varepsilon},
\qquad
\phi_{\varepsilon,\delta}
=\frac{b_\varepsilon-c_\varepsilon+\delta_2K_0}{s_\varepsilon}.
\end{equation}
Then:
\begin{enumerate}
\item $s_\varepsilon\neq0$ and
  $(\theta_{\varepsilon,\delta},\phi_{\varepsilon,\delta})\in\Q^2$;
\item each point $(\theta_{\varepsilon,\delta},\phi_{\varepsilon,\delta})$
  lies in $\cR$ and generates a triple proportional to
  $(a_\varepsilon,b_\varepsilon,c_\varepsilon)$, hence a triangle
  similar to the one generated by $(\theta,\phi)$, with the same
  sides carrying the rational medians;
\item $(\theta_{\varepsilon,\delta},\phi_{\varepsilon,\delta})=(\theta,\phi)$
  for $\varepsilon=\delta=(1,1)$;
\item conversely, every point of $\cR\cap\Q^2$ with the property in
  \textup{(ii)} is one of the points~\eqref{eq:sixteen}.
\end{enumerate}
\end{proposition}

\begin{proof}
(i) The numbers $s_\varepsilon$ are $T_1$, $T_2$, $T_3$ and $-T_4$,
which are nonzero because $H\neq0$; $K_0$ and $L_0$ are polynomials.

(ii) Sign changes of entries do not change $Q$, $K^2$ or $L^2$, so the
triple $(a_\varepsilon,b_\varepsilon,c_\varepsilon)$ satisfies
$Q=H(\theta,\phi)\neq0$ and has doubled medians $\pm K_0$, $\pm L_0$.
The points~\eqref{eq:sixteen} are the pairs of roots
in~\eqref{eq:inverse} for this triple. By
Theorem~\ref{thm:completeness} there is $\rho\neq0$ with
$(A_0,B_0,C_0)(\theta_{\varepsilon,\delta},\phi_{\varepsilon,\delta})
=\rho(a_\varepsilon,b_\varepsilon,c_\varepsilon)$, and then
$H(\theta_{\varepsilon,\delta},\phi_{\varepsilon,\delta})
=\rho^4H(\theta,\phi)>0$.

(iii) This is~\eqref{eq:inverse-K0L0}.

(iv) Let $(\theta',\phi')\in\cR\cap\Q^2$ generate a triple
$(A_0',B_0',C_0')$ with $|A_0'|,|B_0'|,|C_0'|$ proportional to
$|A_0|,|B_0|,|C_0|$. After an overall sign change,
$(A_0',B_0',C_0')=\rho(a_\varepsilon,b_\varepsilon,c_\varepsilon)$
for some $\varepsilon$ and $\rho\neq0$. By~(iii) applied to
$(\theta',\phi')$, the numbers $\theta'$ and $\phi'$ are roots of the
quadratics~\eqref{eq:quadratics} of $(A_0',B_0',C_0')$. These
quadratics are $\rho$ times those of
$(a_\varepsilon,b_\varepsilon,c_\varepsilon)$, whose roots are the
coordinates in~\eqref{eq:sixteen}.
\end{proof}

\begin{remark}\label{rem:involutions}
For $\varepsilon=(1,1)$ the four points are $(\theta,\phi)$,
$(\theta',\phi)$, $(\theta,\phi')$ and $(\theta',\phi')$, where the
product of the roots of~\eqref{eq:quadratics} gives the explicit
involutions
\[
\theta'=\frac{(\phi-1)(\theta+2\phi+1)}{G},\qquad
\phi'=-\frac{(\theta+1)(\phi+2\theta-1)}{G}.
\]
In all our computations the sixteen points were distinct and exactly
one of them lay in $\cF$; for instance,
Euler's triangle (Example~\ref{ex:first-med}) has the point
$(9/16,\,25/52)\in\cF$, and the triangle of
Example~\ref{ex:second-med} has the point $(16/19,\,14/19)\in\cF$.
\end{remark}

% =====================================================================
\section{Local Obstructions to Rational Area}
\label{sec:area}
% =====================================================================

Geometric admissibility guarantees positive area, but not rational
area. Throughout this section assume
\[
(\theta,\phi)\in\cR\cap\Q^2,\qquad \tau\in\Q\setminus\{0\}.
\]
In particular every factor in~\eqref{eq:H} is nonzero. The results
do not use the third-median condition.

\subsection{The rational-area criterion}

Let $S$ be the area of the triangle with sides
$|\tau A_0|,|\tau B_0|,|\tau C_0|$.

\begin{proposition}[Rational-area criterion]
\label{prop:area-criterion}
The generated triangle has rational area if and only if
$H(\theta,\phi)$ is a rational square. This condition does not depend
on $\tau$.
\end{proposition}

\begin{proof}
By Heron's formula, sign invariance and homogeneity,
\[
16S^2=Q(\tau A_0,\tau B_0,\tau C_0)=\tau^4H(\theta,\phi),
\qquad\text{so}\qquad
S=\frac{\tau^2}{4}\sqrt{H(\theta,\phi)}.\qedhere
\]
\end{proof}

A rational square is positive and is a square in $\Q_2$. Positivity
of $H$ is signed admissibility. For the $2$-adic condition, let
$\Zt=\{r/s\in\Q: s\text{ odd}\}$; congruences between elements of
$\Zt$ are understood in this ring. For $x\in\Q^\times$ write
$x=2^{v_2(x)}x^\ast$ with $x^\ast\in\Zt^\times$. If $x$ is a rational
square, then $v_2(x)$ is even and $x^\ast\equiv1\pmod8$, because
$x^\ast$ is then a quotient of odd squares. Every
$u\in\Zt^\times$ satisfies $u^2\equiv1\pmod8$.

\subsection{Valuation obstructions}

Write
\[
\mu=v_2(\theta),\quad \nu=v_2(\phi),\quad
\kappa=\min\{\mu,\nu\},\quad \lambda=\max\{\mu,\nu\}.
\]
We use $v_2(xy)=v_2(x)+v_2(y)$, and $v_2(x_1+\dots+x_k)=v_2(x_i)$
whenever $x_i$ is the unique summand of smallest valuation.

\begin{proposition}[Valuation obstructions]
\label{prop:valuation-obstructions}
In the following three ranges,
\[
v_2(H)=
\begin{cases}
4+\mu+\nu, &\mathrm{(P)}\quad \mu,\nu\ge1,\\
5+7\kappa+\lambda, &\mathrm{(M)}\quad \kappa\le-2,\ \lambda\ge1,\\
5+7\kappa+4\lambda, &\mathrm{(N)}\quad \lambda\le-1,\ \kappa\le\lambda-2.
\end{cases}
\]
Consequently, the generated triangle has irrational area, under every
nonzero rational scaling, in each of the following cases:
\[
\mathrm{(P)}\ \ \mu+\nu\ \text{odd};\qquad
\mathrm{(M)}\ \ \mu+\nu\ \text{even};\qquad
\mathrm{(N)}\ \ \kappa\ \text{even}.
\]
\end{proposition}

\begin{proof}
By~\eqref{eq:H}, $v_2(H)$ is $4$ plus the sum of the valuations of the
ten factors $\phi$, $1-\phi$, $\phi+1$, $G$, $\theta$, $\theta-1$,
$\theta+1$, $\phi-\theta-1$, $\theta+2\phi+1$, $\phi+2\theta-1$.
Table~\ref{tab:valuations} lists these valuations; (M) and (N) are
split according to which parameter has the smaller valuation.

\begin{table}[ht]
\centering
\small
\setlength{\tabcolsep}{4.5pt}
\renewcommand{\arraystretch}{1.15}
\begin{tabular}{lccccc}
\toprule
& (P) & \multicolumn{2}{c}{(M)} & \multicolumn{2}{c}{(N)}\\
\cmidrule(lr){3-4}\cmidrule(lr){5-6}
Factor & $\mu,\nu\ge1$ & $\mu\le-2,\ \nu\ge1$ & $\mu\ge1,\ \nu\le-2$
& $\mu\le\nu-2$ & $\nu\le\mu-2$\\
\midrule
$\phi$            & $\nu$ & $\nu$   & $\nu$   & $\nu$     & $\nu$\\
$1-\phi$          & $0$   & $0$     & $\nu$   & $\nu$     & $\nu$\\
$\phi+1$          & $0$   & $0$     & $\nu$   & $\nu$     & $\nu$\\
$G$               & $0$   & $\mu$   & $\nu$   & $\mu+\nu$ & $\mu+\nu$\\
$\theta$          & $\mu$ & $\mu$   & $\mu$   & $\mu$     & $\mu$\\
$\theta-1$        & $0$   & $\mu$   & $0$     & $\mu$     & $\mu$\\
$\theta+1$        & $0$   & $\mu$   & $0$     & $\mu$     & $\mu$\\
$\phi-\theta-1$   & $0$   & $\mu$   & $\nu$   & $\mu$     & $\nu$\\
$\theta+2\phi+1$  & $0$   & $\mu$   & $\nu+1$ & $\mu$     & $\nu+1$\\
$\phi+2\theta-1$  & $0$   & $\mu+1$ & $\nu$   & $\mu+1$   & $\nu$\\
\midrule
$v_2(H)$ & $4+\mu+\nu$ & $5+7\mu+\nu$ & $5+\mu+7\nu$
& $5+7\mu+4\nu$ & $5+4\mu+7\nu$\\
\bottomrule
\end{tabular}
\caption{Valuations of the factors of $H$; in (N) both valuations are
negative.}
\label{tab:valuations}
\end{table}

Each entry is the valuation of the unique summand of smallest
valuation. For the factors that are linear in $\theta$ and $\phi$ this
is immediate once one notes $v_2(2\theta)=\mu+1$ and
$v_2(2\phi)=\nu+1$. The assumptions $\kappa\le-2$ in (M) and
$|\mu-\nu|\ge2$ in (N) rule out cancellation: for instance,
in (M) with $\mu\le-2$ we have $v_2(2\theta)\le-1<0$, so $2\theta$
cannot cancel the constant in $\phi+2\theta-1$; in (N) with
$\mu\le\nu-2$ we have $\mu+1<\nu$, so $2\theta$ is the unique smallest
summand of $\phi+2\theta-1$. For the mixed factor we use the
appropriate one of
\[
G=(3\phi+1)\theta+(1-\phi)=(3\theta-1)\phi+(\theta+1)
=3\theta\phi+\theta-\phi+1 .
\]
In (P), $v_2((3\phi+1)\theta)=\mu\ge1$ and $1-\phi$ is a unit, so
$v_2(G)=0$. In (M) with $\mu\le-2$, the coefficient $3\phi+1$ is a
unit, so $v_2(G)=\mu$; with $\nu\le-2$, the coefficients $3\theta-1$
and $\theta+1$ are units, so $v_2(G)=\nu$. In (N) both valuations are negative, so the
summand $3\theta\phi$ has valuation $\mu+\nu$, smaller than $\mu$,
$\nu$ and $0$.

Adding the columns gives the stated formulas, and
$5+7\mu+\nu$, $5+\mu+7\nu$ both equal $5+7\kappa+\lambda$, while
$5+7\mu+4\nu$, $5+4\mu+7\nu$ both equal $5+7\kappa+4\lambda$. Modulo
$2$ these are $\mu+\nu$, $1+\mu+\nu$ and $1+\kappa$. In each of the
listed cases $v_2(H)$ is odd, so $H$ is not a rational square, and
Proposition~\ref{prop:area-criterion} completes the proof.
\end{proof}

\begin{remark}
Positive valuation allows arbitrary odd denominators, and negative
valuation describes denominators divisible by a power of $2$, so the
obstructions are not restricted to integral parameters. The cases not
covered are those in which a valuation equals $0$, the mixed case with
$\kappa=-1$, and the case of two negative valuations with
$|\mu-\nu|\le1$. In these cases summands of equal smallest valuation
can cancel; for instance, when $\mu=-1$ and $\nu\ge1$ the summand
$2\theta$ has valuation $0$ and may cancel the constant term of
$\phi+2\theta-1$.
\end{remark}

\subsection{Unit-part obstructions}

Write $\theta=2^\mu u_\theta$ and $\phi=2^\nu u_\phi$ with
$u_\theta,u_\phi\in\Zt^\times$.

\begin{proposition}[Positive-range unit obstruction]
\label{prop:positive-unit}
Suppose $\mu,\nu\ge2$. Rational area requires $\mu+\nu$ to be even and
\[
u_\theta u_\phi\equiv
\begin{cases}
7\pmod8, &\mu=\nu=2\ \text{or}\ \mu,\nu\ge3,\\
3\pmod8, &\text{exactly one of }\mu,\nu\text{ equals }2.
\end{cases}
\]
If either requirement fails, the area is irrational under every
nonzero rational scaling.
\end{proposition}

\begin{proof}
Write $H=16\theta\phi\,U$ with
\[
U=(1-\phi)(\phi+1)\,G\,(\theta-1)(\phi-\theta-1)
(\theta+2\phi+1)(\theta+1)(\phi+2\theta-1).
\]
Since $4\mid\theta$ and $4\mid\phi$ in $\Zt$, we have
$\theta^2\equiv\phi^2\equiv\theta\phi\equiv2\theta\equiv2\phi\equiv0
\pmod8$. Hence, modulo $8$,
\[
\begin{aligned}
(1-\phi)(\phi+1)&\equiv1, &
G&\equiv1+\theta-\phi,\\
(\theta-1)(\phi-\theta-1)&\equiv1-\phi, &
(\theta+2\phi+1)(\theta+1)&\equiv1,\\
\phi+2\theta-1&\equiv\phi-1,
\end{aligned}
\]
and, using $(1-\phi)^2\equiv1$, we get $U\equiv-1-\theta+\phi\pmod8$.
A parameter of valuation exactly $2$ is $\equiv4\pmod8$, and one of
valuation at least $3$ is $\equiv0$. Therefore $U\equiv7$ if
$\mu=\nu=2$ or $\mu,\nu\ge3$, and $U\equiv3$ if exactly one valuation
equals $2$. In particular $U$ is a unit, and
$H=2^{4+\mu+\nu}u_\theta u_\phi U$. A rational square needs an even
exponent and $u_\theta u_\phi U\equiv1\pmod8$; since $3$ and $7$ are
their own inverses modulo $8$, the stated conditions follow.
\end{proof}

\begin{proposition}[Mixed-range unit obstruction]
\label{prop:mixed-unit}
Suppose $\kappa\le-3$ and $\lambda\ge3$. Rational area requires
$\mu+\nu$ to be odd and
\[
u_\theta u_\phi\equiv
\begin{cases}
3\pmod8, &\kappa=-3,\\
7\pmod8, &\kappa\le-4.
\end{cases}
\]
If either requirement fails, the area is irrational under every
nonzero rational scaling.
\end{proposition}

\begin{proof}
First let $\mu\le-3$ and put $\zeta=1/\theta$. Dividing each factor
of~\eqref{eq:H} that involves $\theta$ by $\theta$ (and
$\phi+2\theta-1$ by $2\theta$) gives $H=32\theta^7\phi\,V_\theta$ with
\[
\begin{aligned}
V_\theta={}&(1-\phi^2)\big(1+3\phi+(1-\phi)\zeta\big)(1-\zeta)
\big(-1+(\phi-1)\zeta\big)\\
&\times\big(1+(2\phi+1)\zeta\big)(1+\zeta)
\Big(1+\frac{(\phi-1)\zeta}{2}\Big).
\end{aligned}
\]
Both $\phi$ and $\zeta$ are divisible by $8$ in $\Zt$, and
$v_2(\phi\zeta/2)\ge5$, so $V_\theta\equiv-(1-\zeta/2)\pmod8$.

Next let $\nu\le-3$ and put $\zeta=1/\phi$. In the same way
$H=32\theta\phi^7\,V_\phi$ with
\[
\begin{aligned}
V_\phi={}&(\zeta-1)(1+\zeta)\big(3\theta-1+(\theta+1)\zeta\big)
(\theta^2-1)\\
&\times\big(1-(\theta+1)\zeta\big)
\Big(1+\frac{(\theta+1)\zeta}{2}\Big)\big(1+(2\theta-1)\zeta\big),
\end{aligned}
\]
and $V_\phi\equiv-(1+\zeta/2)\pmod8$.

In both cases $v_2(\zeta/2)=-\kappa-1$, so $\zeta/2\equiv4\pmod8$ if
$\kappa=-3$ and $\zeta/2\equiv0$ if $\kappa\le-4$. Hence the
normalized factor $V$ is a unit with $V\equiv3$ if $\kappa=-3$ and
$V\equiv7$ if $\kappa\le-4$. The valuation of $H$ is
$5+7\kappa+\lambda$, which is even exactly when $\mu+\nu$ is odd, and
its unit part is $u_\theta^7u_\phi V_\theta$ or
$u_\theta u_\phi^7V_\phi$. Since $u^7\equiv u\pmod8$ for units, this is
$\equiv u_\theta u_\phi V\pmod8$, and the condition
$u_\theta u_\phi V\equiv1$ gives the stated residues.
\end{proof}

\begin{remark}[Numerator--denominator form]
Write $\theta=2^\mu m_1/n_1$ and $\phi=2^\nu m_2/n_2$ in lowest terms
with $m_1,m_2,n_1,n_2$ odd and $n_1,n_2>0$. Then
$u_\theta u_\phi\equiv c\pmod8$ is equivalent to
$m_1m_2\equiv c\,n_1n_2\pmod8$, so the unit conditions allow arbitrary
odd denominators.
\end{remark}

\subsection{Using all sixteen parameter points}

\begin{corollary}\label{cor:sixteen-tests}
Let $(\theta,\phi)\in\cR\cap\Q^2$. If one of the sixteen points
of Proposition~\ref{prop:sixteen} satisfies an exclusion condition
of Proposition~\ref{prop:valuation-obstructions},
\ref{prop:positive-unit} or~\ref{prop:mixed-unit}, then the triangle
generated by $(\theta,\phi)$ has irrational area under every nonzero
rational scaling.
\end{corollary}

\begin{proof}
By Proposition~\ref{prop:sixteen}(ii) all sixteen points generate
similar triangles, and similar rational triangles have rational area
simultaneously.
\end{proof}

\begin{example}
Euler's triangle (Example~\ref{ex:first-med}) has $P_1=(11/2,\,2)$
with $(\mu,\nu)=(-1,1)$, which is not covered by
Proposition~\ref{prop:valuation-obstructions}. Its parameter point
$(9/16,\,25/52)\in\cF$ has $(\mu,\nu)=(-4,-2)$, so (N) applies with
$\kappa=-4$ even, and the area is irrational. For
Example~\ref{ex:second-med}, both $P_2=(-6/11,\,16)$ with
$(\mu,\nu)=(1,4)$ and its point $(16/19,\,14/19)\in\cF$ with
$(\mu,\nu)=(4,1)$ are excluded by (P).
\end{example}

\begin{remark}\label{rem:sixteen-numerics}
Using all sixteen points makes the tests considerably stronger. In a
random sample of $2183$ signed admissible points $(\theta,\phi)$ with
numerators of absolute value at most $50$ and denominators at most
$30$, the tests at $(\theta,\phi)$ alone exclude $11.7\%$ of the
points, while the tests at all sixteen points exclude $45.8\%$. As a
consistency check, we applied all tests at the $144$ parameter points
of the nine smallest known Heron triangles with two rational
medians~\cite[Table~1]{Hone2022}; as it must be, none of them is
excluded, and the unit tests of Propositions~\ref{prop:positive-unit}
and~\ref{prop:mixed-unit} apply non-trivially at six of these points.
\end{remark}

\subsection{Explicit infinite families}

\begin{corollary}[Infinite signed admissible families]
\label{cor:infinite-families}
For every integer $n\ge0$, each of the following parameter points is
signed admissible and generates a triangle of irrational area:
\[
\begin{aligned}
&\mathrm{(i)}\ \ \theta=8(2n+1),\ \phi=4; &
&\mathrm{(ii)}\ \ \theta=-\tfrac{1}{4(2n+1)},\ \phi=4;\\
&\mathrm{(iii)}\ \ \theta=4,\ \phi=-\tfrac{1}{4(2n+1)}; &
&\mathrm{(iv)}\ \ \theta=8(8n+1),\ \phi=8;\\
&\mathrm{(v)}\ \ \theta=-\tfrac{1}{8(8n+1)},\ \phi=16; &
&\mathrm{(vi)}\ \ \theta=16,\ \phi=-\tfrac{1}{8(8n+1)};\\
&\mathrm{(vii)}\ \ \theta=\tfrac{1}{16(2n+1)},\ \phi=-\tfrac14 .
\end{aligned}
\]
Families \textup{(i)}--\textup{(iii)} and \textup{(vii)} fail a
valuation condition; families \textup{(iv)}--\textup{(vi)} pass the
relevant parity condition but fail a unit condition.
\end{corollary}

\begin{proof}
\emph{Arithmetic.} The valuations $(\mu,\nu)$ are $(3,2)$ in (i),
$(-2,2)$ in (ii) and $(2,-2)$ in (iii); these are excluded by (P) and
(M) of Proposition~\ref{prop:valuation-obstructions}. In (vii)
$(\mu,\nu)=(-4,-2)$, which is excluded by (N). In (iv)
$(\mu,\nu)=(3,3)$ and $u_\theta u_\phi=8n+1\equiv1\not\equiv7$, so
Proposition~\ref{prop:positive-unit} applies. In (v) and (vi)
$(\mu,\nu)=(-3,4)$ and $(4,-3)$, the sum $\mu+\nu$ is odd, and
$u_\theta u_\phi=-1/(8n+1)\equiv7\not\equiv3$, so
Proposition~\ref{prop:mixed-unit} applies with $\kappa=-3$.

\emph{Geometry.} We count negative factors in~\eqref{eq:H}.
In (i) and (iv), $\theta\ge\phi>1$; exactly $1-\phi$ and
$\phi-\theta-1$ are negative. In (ii) and (v),
$-\frac14\le\theta<0$ and $\phi\in\{4,16\}$; then
$G=(3\phi+1)\theta+1-\phi<0$, and exactly $1-\phi$, $G$, $\theta-1$
and $\theta$ are negative. In (iii), $\theta=4$ and
$-\frac14\le\phi<0$, so $G=5+11\phi\ge\frac94$; in (vi), $\theta=16$
and $-\frac18\le\phi<0$, so $G=17+47\phi\ge\frac{89}{8}$; in both,
exactly $\phi$ and $\phi-\theta-1$ are negative. In (vii),
$0<\theta\le\frac1{16}$ and $\phi=-\frac14$, so $G=\frac{\theta}{4}
+\frac54>0$ and $\theta+2\phi+1=\theta+\frac12>0$; exactly $\phi$,
$\theta-1$, $\phi-\theta-1$ and $\phi+2\theta-1$ are negative. In
every case the number of negative factors is even, so $H>0$.
\end{proof}

These are families of geometrically admissible parameters; we do not
claim that their points satisfy the third-median condition. If such a
point does, its irrational area still prevents a perfect
representative. The tests are necessary rather than sufficient:
passing them does not show that $H$ is a rational square, nor that
the third median is rational.

\begin{remark}[Mixed signs are needed]\label{rem:signed-needed}
The points of families \textup{(i)}--\textup{(vi)} have mixed signs:
they lie outside $\cD\subset(-1,1)^2$, and their entries
$A_0,B_0,C_0$ do not all have the same sign. If one works only with
positive side lengths, the available parameter points of a triangle
are the four points given by~\eqref{eq:inverse} for the positive
triple; these are the points of Proposition~\ref{prop:sixteen} whose
entries have a common sign. For families \textup{(i)}, \textup{(ii)},
\textup{(iv)}, \textup{(v)} and \textup{(vi)} and $0\le n<200$, we
checked that none of these four points satisfies an exclusion
condition of Proposition~\ref{prop:valuation-obstructions},
\ref{prop:positive-unit} or~\ref{prop:mixed-unit}. For instance, the
triangle of family \textup{(ii)} with $n=0$ has the point
$(8/13,\,7/13)\in\cF$ with $(\mu,\nu)=(3,0)$, which is not covered by
these propositions, whereas its mixed-sign point $(-1/4,\,4)$ is
excluded by \textup{(M)}. For these families, the simple valuation
criteria apply only at mixed-sign parameter points.
\end{remark}

% =====================================================================
\section{Conclusion}
\label{sec:conclusion}
% =====================================================================

Buchholz's parametrization separates three requirements: geometric
admissibility, rationality of the third median, and rationality of
the area. The factored Heron polynomial describes the first
completely: signed admissibility is the condition $H>0$, and direct
admissibility is the union of three quadrilaterals containing the
Buchholz--Rathbun region. Every triangle with two rational medians
corresponds to sixteen explicit parameter points in the signed region.

The same factorization gives $2$-adic obstructions to rational area:
parity conditions in three valuation ranges and unit conditions
modulo $8$ in the positive and deeper mixed ranges. Applied at all
sixteen parameter points of a triangle, they exclude infinite
families of admissible parameters, including non-integral ones; for
several of these families the exclusion is visible only at
mixed-sign parameter points. These
conditions do not classify the rational-area parameters, and they do
not resolve the existence of perfect triangles. Natural next steps are
the valuation ranges not covered here, and analogous conditions at odd
primes.

\subsection*{Acknowledgements}
% TODO: add thanks to your supervisor here.
% TODO: before posting, make sure that the AI statement below matches
% the final version and that you have checked every proof yourself.
The author used Claude (Anthropic), an AI assistant, to check
computations, to edit the exposition, and to suggest
Lemma~\ref{lem:cross}, Theorem~\ref{thm:completeness}, the description
of $\cD$ in Proposition~\ref{prop:direct-region},
Proposition~\ref{prop:sixteen}, case (N) of
Proposition~\ref{prop:valuation-obstructions}, family (vii) of
Corollary~\ref{cor:infinite-families}, and the computations in
Remarks~\ref{rem:sixteen-numerics} and~\ref{rem:signed-needed}. All
statements and proofs were checked by the author, who takes full
responsibility for the content.

% =====================================================================
\appendix
\section{Buchholz's derivation}
\label{app:derivation}
% =====================================================================

For completeness we sketch how the parametrization arises; nothing in
the main text depends on this appendix. Let $a,b,c,K,L\in\Q$
satisfy~\eqref{eq:doubled-medians}. The first relation is
$K^2+a^2=2(b^2+c^2)$, so in $\Q(\sqrt2)$
\[
(K-b\sqrt2)(K+b\sqrt2)=-(a-c\sqrt2)(a+c\sqrt2).
\]
Multiplying by $a-c\sqrt2$ and using $a^2-2c^2=-(K^2-2b^2)$ gives
\[
\frac{K+b\sqrt2}{a+c\sqrt2}=\frac{r_1+s_1\sqrt2}{t_1},
\qquad
r_1=2bc-aK,\quad s_1=cK-ab,\quad t_1=K^2-2b^2 .
\]
Taking norms, $(r_1^2-2s_1^2)/t_1^2=-1$, that is
$r_1^2+t_1^2=2s_1^2$. After clearing denominators $r_1\equiv t_1
\pmod 2$, and $R=(r_1+t_1)/2$, $T=(r_1-t_1)/2$ satisfy $R^2+T^2=s_1^2$.
The parametrization of Pythagorean triples gives, up to a common factor
and a choice of signs,
\begin{equation}\label{eq:pyth}
r_i=2p_iq_i+p_i^2-q_i^2,\qquad
s_i=p_i^2+q_i^2,\qquad
t_i=2p_iq_i-p_i^2+q_i^2
\end{equation}
with $i=1$. In the same way $L^2+b^2=2(a^2+c^2)$ leads
to~\eqref{eq:pyth} with $i=2$, and we obtain
\[
K+b\sqrt2=\frac{r_1+s_1\sqrt2}{t_1}\,(a+c\sqrt2),\qquad
L+c\sqrt2=\frac{r_2+s_2\sqrt2}{t_2}\,(b+a\sqrt2).
\]
Comparing the coefficients of $\sqrt2$ gives the linear relations
\[
r_1c+s_1a=t_1b,\qquad r_2a+s_2b=t_2c,
\]
whose solution is, up to a common factor,
\[
(a,b,c)=\big(t_1t_2-r_1s_2,\ s_1t_2+r_1r_2,\ t_1r_2+s_1s_2\big)
=2q_1^2q_2^2\,\big(A_0,B_0,C_0\big)(\theta,\phi),
\]
with $\phi=p_1/q_1$ and $\theta=p_2/q_2$. By
Theorem~\ref{thm:completeness}, the sign choices made
in~\eqref{eq:pyth} entail no loss of generality.

% =====================================================================

\end{document}